# MuLoG, or How to apply Gaussian denoisers to multi-channel SAR speckle reduction?

Charles-Alban Deledalle, Loïc Denis, Sonia Tabti, Florence Tupin,

*Abstract*—Speckle reduction is a longstanding topic in synthetic aperture radar (SAR) imaging. Since most current and planned SAR imaging satellites operate in polarimetric, interferometric or tomographic modes, SAR images are multi-channel and speckle reduction techniques must jointly process all channels to recover polarimetric and interferometric information. The distinctive nature of SAR signal (complex-valued, corrupted by multiplicative fluctuations) calls for the development of specialized methods for speckle reduction. Image denoising is a very active topic in image processing with a wide variety of approaches and many denoising algorithms available, almost always designed for additive Gaussian noise suppression. This paper proposes a general scheme, called MuLoG (MUlti-channel LOgarithm with Gaussian denoising), to include such Gaussian denoisers within a multi-channel SAR speckle reduction technique. A new family of speckle reduction algorithms can thus be obtained, benefiting from the ongoing progress in Gaussian denoising, and offering several speckle reduction results often displaying method-specific artifacts that can be dismissed by comparison between results.

## I. INTRODUCTION

Synthetic aperture radar (SAR) imaging is a widely used technique for earth observation. It offers complementary information to the more common optical imaging. Among its distinctive features, one may cite its all-weather and day-and-night imaging capabilities [50], the interferometric configurations that give access to 3-D reconstructions for digital elevation models [55] and displacement estimation [49], [25] or the polarimetric and tomographic modes that give access to estimates of the biomass in forested areas [39].

Coherent combination of several radar echoes within each resolution cell results in interferences and the well-known *speckle phenomenon* [31]. Due to speckle, regions with homogeneous radar properties display strong fluctuations in SAR images. Direct estimation of the reflectivity, the interferometric phase or polarimetric properties is unusable given its prohibitively large variance. Speckle reduction is thus a longstanding topic in SAR imagery.

Since speckle statistics departs from the additive Gaussian noise model widely used for processing optical images, a whole line of denoising methods have been designed specifically for speckle reduction. Restoration of intensity images (i.e., single-channel SAR images) remains the first and most studied case. Many different schemes have been proposed in this context, see the recent reviews [3] and [15]. Among the different possible strategies, *selection-based methods* use various criteria to identify homogeneous collections of neighboring pixels. Pixels may be selected by locally choosing the best neighborhood within a fixed set of oriented windows [41]. Another approach consists of selecting connected pixels by region growing [64]. Patch-comparison has been shown to provide a robust way to get a (weighted) selection of similar pixels [19], [15]. *Variational methods* formulate the estimation problem as an optimization problem depending on the whole image. The objective function that is optimized is composed of two terms: a data-fitting term and a regularization term. Due to the non-Gaussian distribution of the intensity of SAR images, the data-fitting term differs from the usual sum of squared differences that arises from a Gaussian assumption. Several regularization terms have been considered in the literature, e.g., total variation (TV) [57], [59], [22], [60], [7], curvelets [23], Gaussian mixture models [63]. Yet another strategy is to transform the data so the noise becomes additive and approximately Gaussian and stationary by using a *homomorphic transform*. As further discussed in Sec. II, the log-transformed intensity is approximately Gaussian distributed and can thus be restored using a Gaussian denoiser, e.g., based on wavelet thresholding [68], [1], [6] or on patch redundancy [48].

Although many applications are based on the intensity provided by a single SAR image (corresponding to the square of the modulus of the backscattered electromagnetic field), the SAR data can carry much more information. Indeed, the phase of a pixel carries information about the acquisition geometry. Combining two SAR acquisitions with a slightly different angle allows retrieving information on the elevation of the points and corresponds to the so-called interferometric mode. Processing this information can be done by considering the vector of the two complex SAR signals and their associated empirical covariance matrix. Another information is provided by the polarization of the emitted and backscattered wave. By emitting and receiving the signal with different polarizations (for instance vertical and horizontal polarizations), a vector of 4 complex values is recorded for each pixel. This vector gives valuable information on the backscattering mechanisms happening inside the pixel (double bounce scattering, volume scattering, etc.). This polarimetric mode can also be combined with the interferometric one, giving a D-dimensional complex vector for each pixel and allowing to retrieve the height of

C. Deledalle is with IMB, CNRS, Univ Bordeaux, Bordeaux INP, F-33405 Talence, France, e-mail: charles-alban.deledalle@math.u-bordeaux.fr.

L. Denis is with Univ Lyon, UJM-Saint-Etienne, CNRS, Institut d'Optique Graduate School, Laboratoire Hubert Curien UMR 5516, F-42023, Saint-Etienne, France, e-mail: loic.denis@univ-st-etienne.fr

S. Tabti was with LTCI, Telecom ParisTech, Université Paris-Saclay, 75013, Paris, France, email: sonia.tabti@telecom-paristech.fr

F. Tupin is with LTCI, Telecom ParisTech, Université Paris-Saclay, 75013, Paris, France, email: florence.tupin@telecom-paristech.fr

different scatterers inside the pixel. Finally, the tomographic mode exploits a set of interferometric acquisitions to recover the elevation profile inside the pixel. These multi-channel SAR data are of high interest for a wide range of applications for Earth observation. But as for the intensity data, the speckle phenomenon induces a high variability of the physical parameters (interferometric phase, coherence, inter-channel cross-correlations, etc.).

Not all the aforementioned speckle reduction methods generalize well to multi-channel SAR images. *Selection-based methods* require extending the homogeneity or similarity criteria to multi-variate data. This can be done by using only part of the information, e.g., the span [64] or the scattering properties [42], or by exploiting the whole covariance matrices [10], [20]. Once relevant pixels have been selected, estimation of polarimetric/interferometric information is straightforward, e.g., by using a linear mean square error approach [40], [20] or a weighted maximum likelihood estimator [17], [18]. Extension of *variational methods* to multi-channel SAR data is more challenging. A simple remedy is to consider only restoring the diagonal elements [44] at a price of information loss about cross-channel correlations. Direct formulation of the objective function on the full covariance matrices raises several problems: (i) computational complexity due to the nonconvexity of the data-fitting term; (ii) difficulty to express regularity properties of the complex-valued terms of the covariance matrices; (iii) non-stationary variance of speckle that leads to over/under-smoothing in some areas. These difficulties explain that very few works were conducted in this direction, with the exception of recent works on multi-channel TV regularizations [52], [53], [54]. Finally, while the *homomorphic transform* approach is well understood on intensity images, no variance stabilization transform is known for multi-channel images.

A generic methodology to apply denoising methods from the "additive Gaussian noise world" to multi-channel SAR data is definitely lacking. This paper attempts to fill this gap. The contributions of this paper are the following:

1) to provide a generic method, called MuLoG, to embed a Gaussian denoiser in a speckle reduction filter;
2) to apply as well to single-channel or to multi-channel SAR images;
3) to produce estimates of the complex covariance matrices that capture all polarimetric and/or interferometric information,
4) to better account for speckle statistics than other methods based on a homomorphic transform[1];
5) to require no parameter tuning other than possibly within the Gaussian denoiser.

We introduce our generic methodology by first considering (single-channel) intensity images (Sec. II), then the extension to multi-channel SAR images (Sec. III). We discuss implementation issues (Sec. IV) before illustrating the proposed methodology with several Gaussian denoisers and different types of multi-channel SAR images (Sec. V).

---

[1] note that these other methods are only applicable to intensity SAR images (i.e., single-channel)

## II. Speckle reduction for SAR intensity images

### A. Statistics of univariate SAR images

*a) Intensity:* Univariate SAR images are by nature complex-valued and only the modulus (*a.k.a.*, the amplitude) is informative. The square of the modulus (*a.k.a.*, the intensity) is nevertheless easier to manipulate and, according to Goodman's model [31], it follows a gamma distribution $\mathcal{G}(R; L)$ with a probability density given by

$$p_I(I|R) = \frac{L^L I^{L-1}}{\Gamma(L) R^L} \exp\left(-L\frac{I}{R}\right) , \quad (1)$$

where $I \in \mathbb{R}_+$ is the observed intensity, $R \in \mathbb{R}_+$ is the underlying reflectivity (related to the radar cross-section), $L > 0$ is the number of looks, and $\Gamma$ the gamma function. Note that we denote by $\mathbb{R}_+$ the set of positive real values. The intensity $I$ can be decomposed as a product of the reflectivity $R$ and of a speckle component $S$ distributed under a standard gamma distribution ($S \sim \mathcal{G}(1; L)$):

$$I = R \times S, \quad \mathbb{E}[I] = R \quad \text{and} \quad \text{Var}[I] = \frac{R^2}{L} . \quad (2)$$

As the variance depends on the expectation, fluctuations are said to be signal dependent. The top graph in Fig.2(a) illustrates how a simple rectangle signal (solid curve) gets corrupted by speckle: the gray area shows values between the first and third quartiles, the dots represent the expectation and the dashed line a single noisy realization. The signal-dependent nature of the fluctuations can be observed: the difference between first and third quartiles is larger when the underlying signal values are high. Last but not least, the gamma distribution has a heavier right-tail characterizing the typical bright outliers observed in SAR intensity images.

*b) Logarithm:* The log-transform $y = \log I \in \mathbb{R}$ is often employed to convert multiplicative fluctuations to additive ones. From the gamma distribution (1), $y$ follows the Fisher-Tippett distribution defined as

$$p_y(y|x) = \frac{L^L}{\Gamma(L)} e^{L(y-x)} \exp\left(-Le^{y-x}\right) , \quad (3)$$

where $x = \log R \in \mathbb{R}$. The Fisher-Tippett distribution, denoted by $\mathcal{FT}(x; L)$, models additive corruptions as [69]

$$y = x + s , \quad (4)$$
$$\mathbb{E}[y] = x - \log L + \Psi(L) , \quad (5)$$
$$\text{and} \quad \text{Var}[y] = \Psi(1, L) , \quad (6)$$

where $s \sim \mathcal{FT}(0; L)$, $\Psi(\cdot)$ is the digamma function and $\Psi(\cdot, L)$ is the polygamma function of order $L$. The log transform stabilizes the variance, *i.e.*, the fluctuations are made signal independent. Equation (5) shows that the noise has a non-zero mean. If noise is assumed zero-mean during a speckle-reduction step performed on log-transformed data, a subsequent debiasing step is then necessary.

The bottom graph in Fig.2(a) displays the intensities in decibels (*i.e.*, $10 \cdot \log_{10}$ of the intensity). By comparing with the top curve, variance stabilization can be observed, as well as the bias between the expectation of log-transformed intensities (black dots) and ground truth signal (solid curve).

## B. Homomorphic approach

We now consider $I$ and $R$ as images, such that $I \in \mathbb{R}_+^n$, and $R \in \mathbb{R}_+^n$, where $n$ is the number of pixels. Let $y \in \mathbb{R}^n$ and $x \in \mathbb{R}^n$ be the entry-wise logarithm of $I$ and $R$ respectively, *i.e.*, such that $y_k = \log I_k$ and $x_k = \log R_k$. A classical approach to estimate $R$ is thus to approach the Fisher-Tippett distribution by a non-centered Gaussian distribution, which leads to the following estimation procedure

$$\hat{x} = f_{\Psi(1,L)}(y) + \underbrace{(\log L - \Psi(L))\mathbf{1}_n}_{\text{debiasing}} , \qquad (7)$$

where $f_{\sigma^2} : \mathbb{R}^n \to \mathbb{R}^n$ is a denoiser for images contaminated by zero-mean additive white Gaussian noise $\mathcal{N}(0; \sigma^2 \mathrm{Id})$. Typically, $f$ can be a regularized least-square solver expressed as

$$f_{\sigma^2}(y) \in \underset{x \in \mathbb{R}^n}{\operatorname{argmin}} \frac{1}{2\sigma^2}\|y - x\|^2 + \mathcal{R}(x) , \qquad (8)$$

where, within the Maximum A Posteriori (MAP) framework, $\mathcal{R}(x) = -\log p_x(x)$ is a prior term enforcing some regularity on the solution. Finally, the estimate $\hat{R}$ is defined as $\hat{R}_k = \exp \hat{x}_k$ for all pixel index $k$. This process is summarized on Fig.1, top row. We illustrate the restored 1D signals obtained by applying L2+TV minimization [58] within a homomorphic procedure on Fig.2(b).

The homomorphic approach has been extensively used, *e.g.*, for wavelet prior in [68], for patch-based priors as in KSVD [26] or for non-local filtering [48]. While for large values of $L$ the Fisher-Tippett distribution approaches the Gaussian distribution, the two distributions differ significantly for low values of $L$. In particular, the Fisher-Tippett is asymmetrical with a heavier left-tail. Hence, this approximation often leads to remaining dark stains on the resulting images (see also Fig. 6(d)).

## C. Variational approach

Rather than applying the homomorphic approach in order to reduce the original problem to the Gaussian case, an alternative proposed in [4], [22] is to consider directly the gamma distribution of $I \in \mathbb{R}_+^n$ leading to a MAP solver of the form (Fig.1 second row):

$$\hat{R} \in \underset{R \in \mathbb{R}_+^n}{\operatorname{argmin}} -\log p_I(I|R) + \mathcal{R}(R) , \qquad (9)$$

$$\text{where} \quad -\log p_I(I|R) = L\sum_{k=1}^n \log R_k + \frac{I_k}{R_k} + \text{Cst}.$$

Cst. denoting a constant for the optimization problem. Not only is this minimization constrained to positive-valued images, but the objective function is also nonconvex. This nonconvexity enforces some robustness to the bright outliers originating from the right tail of the gamma distribution. As a consequence, when standard iterative solvers are used (*e.g.*, gradient descents [4], or the forward-backward algorithm [11]), the solution depends on the initialization and the internal parameters of the solver, even if $\mathcal{R}$ is chosen convex. Convexification by replacing the original objective function by its convex hull [36] simplifies the minimization at the cost of a loss of accuracy of the statistical model. In particular, the convex hull does not capture the right tail of the gamma distribution, leading to remaining bright outliers on the resulting images. For some Markov prior regularization, namely convex pairwise regularizations, the global optimum can be obtained in finite time [35]. This optimization method has been applied in the 1D illustration of Fig.2(c). It can be observed that, since large fluctuations in bright areas are more strongly penalized than the small fluctuations in low-level areas, speckle is reduced predominantly in the brightest areas. In the restored signal, the remaining noise variance is made uniform, i.e., signal-independent, in linear scale. Thus, small details with an identical signal to noise ratio will be more likely suppressed if placed in areas of large average value. This phenomenon is also observed on images: speckle noise is reduced more strongly in brighter areas [22], [60], and bright punctual targets are found to be spread out.

## D. Variational approach on log-transformed data

Another alternative proposed in [59], [60], [7], [23] consists of formulating the estimation problem in the log domain. Rather than approximating the likelihood of log-transformed reflectivities by a Gaussian distribution (as done by homomorphic approaches discussed in Sec. II-B), the Fisher-Tippett distribution (3) is considered. The regularization is also expressed on the log-transformed reflectivities $x \in \mathbb{R}^n$, leading to a MAP optimization problem of the form:

$$\hat{x} \in \underset{\hat{x} \in \mathbb{R}^n}{\operatorname{argmin}} -\log p_y(y|x) + \mathcal{R}(x) , \qquad (10)$$

$$\text{where} \quad -\log p_y(y|x) = L\sum_{k=1}^n x_k + e^{y_k - x_k} + \text{Cst}.$$

The final estimate is defined as $\hat{R}_k = \exp \hat{x}_k$ for all pixel index $k$. Note that even though the solutions of problems (9) and (10) are different, their definitions only differ in terms of the prior regularization since $\log p_y(y|x) = \log p_I(\exp(y)|\exp(x)) + \text{Cst}$.. Nevertheless, this change of variable leads to several advantages compared to (9). First, the optimization is unconstrained as $x$ can be any vector of $\mathbb{R}^n$. More importantly, the data fidelity becomes convex. As a consequence, if $\mathcal{R}$ is also chosen convex, solutions depend neither on the initialization nor the choice of the internal parameters of the solver. Such internal parameters only affect the speed of convergence and thus the number of iterations to perform in practice.

The multiplicative image denoising by augmented Lagrangian (MIDAL) algorithm [7] considers a convex regularization, TV, and minimizes (10) using the alternating direction method of multipliers (ADMM) algorithm [30], [28], [24] that repeats, for an internal parameter $\beta > 0$, the updates

$$\hat{z} \leftarrow \underset{z}{\operatorname{argmin}} \frac{\beta}{2}\|z - \hat{x} + \hat{d}\|^2 + \mathcal{R}(z) , \qquad (11)$$

$$\hat{d} \leftarrow \hat{d} + \hat{z} - \hat{x} , \qquad (12)$$

$$\hat{x} \leftarrow \underset{x}{\operatorname{argmin}} \frac{\beta}{2}\|x - \hat{z} - \hat{d}\|^2 - \log p_y(y|x) . \qquad (13)$$

Clearly, the minimization for $z$ in (11) depends on the choice of $\mathcal{R}$ and can be solved by a dedicated regularized least square solver (corresponding to a Gaussian denoiser). This is schematized on Fig.1, third row. The minimization for $x$ in (13) amounts to solving $n$ separable convex problems of the form

$$\underset{x_k}{\mathrm{argmin}} \frac{\beta}{2}(x_k - a_k)^2 + L\left(x_k + e^{y_k - x_k}\right) ,\qquad (14)$$

where $a_k = \hat{z}_k + \hat{d}_k$. The explicit solution is given by the Lambert W functions [7], or can be computed more efficiently with a few iterations of Newton's method as

$$\hat{x}_k \leftarrow \hat{x}_k - \frac{\beta(\hat{x}_k - a_k) + L(1 - e^{y_k - \hat{x}_k})}{\beta + Le^{y_k - \hat{x}_k}} .\qquad (15)$$

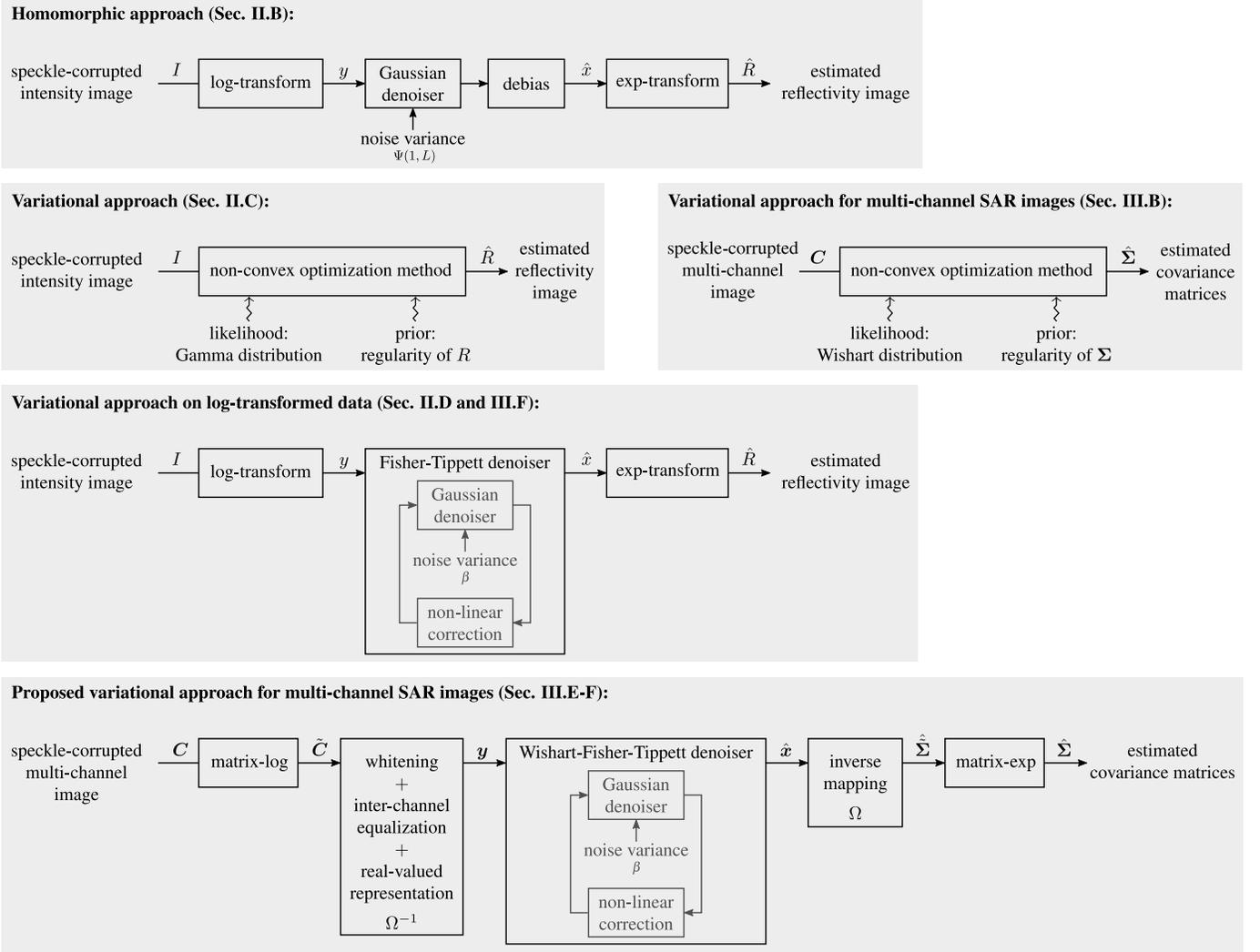

Fig. 1. Several approaches for speckle reduction discussed in this paper. The proposed scheme generalizes the variational approaches based on a log-transform to multi-channel SAR images. Similarly to homomorphic or variational approaches on log-transformed data, it embeds a Gaussian denoising step.

Using ten iterations is usually enough to offer good performances within a reasonable computational time, see, *e.g.*, [63].

As already mentioned, when $\mathcal{R}$ is convex, the parameter $\beta$ only acts on the speed of convergence. Interestingly, as the noise variance is independent of the signal, the convergence is in practice uniform meaning that for a finite number of iterations the same amount of smoothing will be performed both in dark and bright regions. Small details with identical signal to noise ratio will be identically smoothed whatever the average value of the area.

When $\mathcal{R}$ is nonconvex but satisfies some weak conditions and $\beta$ is chosen large enough, ADMM still converges to a local minimum [34]. In this case, the solution depends on both the initialization and the choice of $\beta$. We observe that choosing $\beta$ as $(1 + 2/L)\mathrm{Var}[y]^{-1}$, where $\mathrm{Var}[y] = \Psi(1, L)$, provides a similar smoothing whatever the number of iterations and the initial number of looks $L$.

This shows that taking the logarithm not only makes the data fidelity term convex, but also ensures a uniform speed of convergence of a non-constrained optimization problem. These are two key practical ingredients that challenge nonconvex strategies that directly deal with gamma distributed values.

We illustrate on Fig.2(d) that minimization of the total variation of log-transformed intensities produces a smaller loss of contrast of bright structures and reduces speckle equally in all regions.

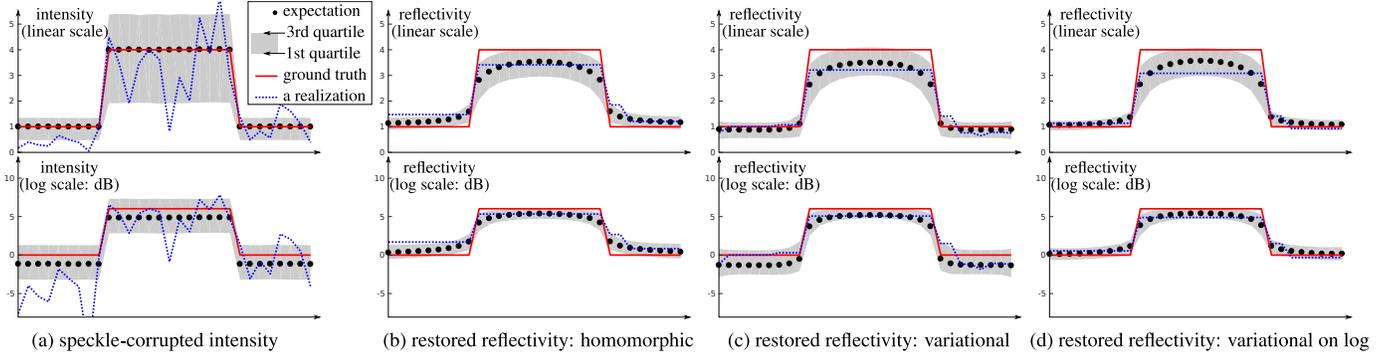

Fig. 2. Illustration of speckle reduction by 1D total variation minimization: (a) speckle-corrupted rectangle signal ($L = 2$), (b) restored reflectivity using the homomorphic approach and "L2 + TV" minimization in log-domain, (c) restored reflectivity using a variational approach "gamma + TV" in the linear domain, (d) restored reflectivity using a variational approach "Fisher-Tippett + TV" in the log domain. Expectation and quartiles are computed from 50 000 noisy realizations. Restorations computed by exact discrete optimization [35] with 2 000 quantization levels.

## III. EXTENSION TO MULTI-VARIATE SAR IMAGERY

### A. Statistics of covariance matrices

As mentioned in the introduction, multi-variate complex SAR images carry much more information than univariate SAR images since inter-channel cross-correlations capture geophysical features (e.g., height is related to the interferometric phase, geometrical configuration with the polarimetric properties). All relevant information of these complex-valued images can be gathered together at each pixel by forming a $D \times D$ complex covariance matrix $\boldsymbol{C}$, with $D = 2$ in single-baseline SAR interferometry, $D = 3$ in SAR polarimetry, $D = 6$ in PolInSAR (polarimetry + interferometry), and even larger values of $D$ for multi-baseline interferometry or SAR tomography [50]. Let us denote by $\vec{k}_i$ the vector of $D$ complex values recorded for a pixel $i$, then the empirical covariance matrix $\boldsymbol{C}$ is given by:

$$\boldsymbol{C} = \frac{1}{L} \sum_{t=1}^{L} \vec{k}_t . \vec{k}_t^*$$

where $^*$ denotes the Hermitian transpose. Note that the $L$ vectors $\vec{k}_t$ are often extracted from a small window centered around the pixel of interest leading to a loss of resolution. If spatial resolution must be preserved, a single-look setting should be used, i.e., $L = 1$, resulting to a matrix $\boldsymbol{C}$ with rank 1. Goodman's model [31] describes speckle in this matrix as being circular complex Wishart distributed, for $L \geq D$:

$$p_{\boldsymbol{C}}(\boldsymbol{C}|\boldsymbol{\Sigma}) = \frac{L^{LD}|\boldsymbol{C}|^{L-D}}{\Gamma_D(L)|\boldsymbol{\Sigma}|^L} \exp\left(-L\operatorname{tr}(\boldsymbol{\Sigma}^{-1}\boldsymbol{C})\right) , \quad (16)$$

where $\boldsymbol{\Sigma}$ is the underlying covariance matrix encoding reflectivities and complex correlations, and $L \geq D$ is the number of looks. Both $\boldsymbol{C}$ and $\boldsymbol{\Sigma}$ belong to the open cone of complex Hermitian positive definite matrices. While fluctuations in univariate SAR images are multiplicative, the Wishart distribution, denoted by $\mathcal{W}(\boldsymbol{\Sigma}; L)$ models fluctuations in multivariate SAR images as

$$\boldsymbol{C} = \boldsymbol{\Sigma}^{1/2}\boldsymbol{S}\boldsymbol{\Sigma}^{1/2} , \quad (17)$$

$$\mathbb{E}[\boldsymbol{C}] = \boldsymbol{\Sigma} , \quad (18)$$

$$\operatorname{Var}[\boldsymbol{C}_{ij}] = \frac{1}{L}\boldsymbol{\Sigma}_{ii}\boldsymbol{\Sigma}_{jj} , \quad (19)$$

where $\boldsymbol{S} \sim \mathcal{W}(\operatorname{Id}; L)$ (see, e.g., [47])[2]. Note that the variance for off-diagonal elements does not depend on $\boldsymbol{\Sigma}_{ij}$ but on $\boldsymbol{\Sigma}_{ii}$ and $\boldsymbol{\Sigma}_{jj}$, which indicates that the fluctuations are not only intra-channel signal dependent but inter-channel signal dependent. Interestingly, according to [32], [47], we have the following relations regarding the determinant and the trace

$$|\boldsymbol{C}| = |\boldsymbol{\Sigma}||\boldsymbol{S}| , \quad (20)$$

$$\operatorname{tr}\boldsymbol{C} = \operatorname{tr}(\boldsymbol{\Sigma}\boldsymbol{S}) , \quad (21)$$

$$\mathbb{E}[\operatorname{tr}\boldsymbol{C}] = \operatorname{tr}\boldsymbol{\Sigma} , \quad (22)$$

$$\operatorname{Var}[\operatorname{tr}\boldsymbol{C}] = \frac{1}{L}\operatorname{tr}\boldsymbol{\Sigma}^2 , \quad (23)$$

and $\mathbb{E}[\operatorname{tr}(\boldsymbol{C}^2 - \boldsymbol{\Sigma}^2)] = \frac{1}{L}(\operatorname{tr}\boldsymbol{\Sigma})^2 . \quad (24)$

### B. Limit of a direct variational approach

As for the univariate case, a variational approach considering the Wishart distribution of $\boldsymbol{C}$ can be expressed as

$$\hat{\boldsymbol{\Sigma}} \in \underset{\boldsymbol{\Sigma} \succ^H 0}{\operatorname{argmin}} \ -\log p_{\boldsymbol{\Sigma}}(\boldsymbol{\Sigma}|\boldsymbol{C}) + \mathcal{R}(\boldsymbol{\Sigma}) , \quad (25)$$

where $-\log p_{\boldsymbol{\Sigma}}(\boldsymbol{\Sigma}|\boldsymbol{C}) = L \sum_{k=1}^{n} \log|\boldsymbol{\Sigma}_k| + \operatorname{tr}(\boldsymbol{\Sigma}_k^{-1}\boldsymbol{C}_k) + \operatorname{Cst}.$

where $\boldsymbol{\Sigma} \succ^H 0$ reads as "$\boldsymbol{\Sigma}$ is Hermitian positive definite". Note that (25) boils down to (9) for $D = 1$. While estimating directly the reflectivity in a variational framework in the case of univariate data requires to enforce a positivity constraint, optimizing for $\boldsymbol{\Sigma}$ requires optimizing on the open cone of complex Hermitian positive definite matrices, which is much more challenging. Moreover, as noted for the univariate case, the neg-log-likelihood associated to the Wishart distribution is highly nonconvex, so that finding a good quality local optimum is very difficult. This approach is summarized on the second row of Fig.1.

To circumvent the difficulty arising from the nonconvexity of the neg-log-likelihood, Nie et al. approximated it by its convex envelope [52], [53]. In their first algorithm, WisTV [52], $\mathcal{R}$ was chosen as a matricial total-variation regularizer,

[2] for complex random variables $\operatorname{Var}[\boldsymbol{C}_{ij}] = \mathbb{E}[|\boldsymbol{C}_{ij}|^2] - |\mathbb{E}[\boldsymbol{C}_{ij}]|^2$

$$\mathcal{K} : \begin{pmatrix} \alpha^1 \\ \vdots \\ \alpha^{D^2} \end{pmatrix} \mapsto \begin{pmatrix} \alpha^1 & (\alpha^{D+1} + j\alpha^{D+2})/\sqrt{2} & \cdots & (\alpha^{D^2-1} + j\alpha^{D^2})/\sqrt{2} \\ (\alpha^{D+1} - j\alpha^{D+2})/\sqrt{2} & \alpha^2 & & \\ \vdots & & \ddots & \vdots \\ (\alpha^{D^2-1} - j\alpha^{D^2})/\sqrt{2} & & \cdots & \alpha^D \end{pmatrix} \tag{26}$$

and next in [53] it was replaced by a matricial non-local total-variation regularizer [29] to better preserve textures. Replacing the neg-log-likelihood by its convex envelope is a rather crude approximation since it leads to underestimating the right tail of the distribution, and thus prevents from being robust against bright outliers.

The noise components in $\boldsymbol{C}$ are intensively signal dependent, so that speckle suppression is not as strong in all regions. Last but not least, as for the univariate case, it has been recently observed that the performance of WisTV can drop significantly when the dynamic range of the input SAR image becomes very high [54]. In this latter article, the authors solve this issue by normalizing each input matrix $\boldsymbol{C}_k$ by a factor $\tau_k$ defined as a non-linear sigmoid function of their span $\frac{1}{3}\operatorname{tr}\boldsymbol{C}_k$. The output $\tilde{\boldsymbol{\Sigma}}_k$ are then rescaled back to their original range by multiplication by $\tau_k^{-1}$. In addition, they reformulate (25) into an unconstrained optimization problem by writing $\boldsymbol{C}_k = \boldsymbol{\Pi}_k \boldsymbol{\Pi}_k^* + \epsilon \mathrm{Id}$, for a fixed $\epsilon > 0$, and then optimizing for $\boldsymbol{\Pi}_k \in \mathbb{R}^{D \times D}$. The factorization for $\boldsymbol{\Pi}_k$ being non unique, the resulting energy will necessarily present even more local minima. To cope with these different issues we have chosen a different path. We suggest mimicking the univariate case by making use of the matrix logarithm, thereby extending MIDAL approach [7] to multi-channel SAR images.

### C. The Wishart-Fisher-Tippett distribution

Our objective is to generalize the use of the log transform in univariate SAR images to multi-variate SAR images. To that end, we resort to the matrix logarithm defined[3,4] as

$$\boldsymbol{\Sigma} \mapsto \tilde{\boldsymbol{\Sigma}} = \log\boldsymbol{\Sigma} = \mathbf{E}\operatorname{diag}(\tilde{\boldsymbol{\Lambda}})\mathbf{E}^{-1} \quad \text{where} \quad \tilde{\boldsymbol{\Lambda}}_i = \log\boldsymbol{\Lambda}_i , \tag{27}$$

$\mathbf{E} \in \mathbb{C}^{D \times D}$ is the matrix whose column vectors are eigenvectors (with unit norm) of $\boldsymbol{\Sigma}$, $\boldsymbol{\Lambda} \in \mathbb{R}_+^D$ is the vector of corresponding eigenvalues, such that $\boldsymbol{\Sigma} = \mathbf{E}\operatorname{diag}(\boldsymbol{\Lambda})\mathbf{E}^{-1}$, and $\tilde{\boldsymbol{\Lambda}} \in \mathbb{R}^D$. Its inverse transform is the matrix exponential defined similarly as

$$\tilde{\boldsymbol{\Sigma}} \mapsto \boldsymbol{\Sigma} = e^{\tilde{\boldsymbol{\Sigma}}} = \mathbf{E}\operatorname{diag}(\boldsymbol{\Lambda})\mathbf{E}^{-1} \quad \text{where} \quad \boldsymbol{\Lambda}_i = \exp\tilde{\boldsymbol{\Lambda}}_i . \tag{28}$$

While $\boldsymbol{\Sigma}$ lies in the open cone of complex Hermitian positive definite matrices, $\tilde{\boldsymbol{\Sigma}}$ lies in the vector space of complex Hermitian matrices which is isomorphic to $\mathbb{R}^{D^2}$. The change of variables $\tilde{\boldsymbol{C}} = \log\boldsymbol{C}$ and $\tilde{\boldsymbol{\Sigma}} = \log\boldsymbol{\Sigma}$ leads to the distribution of log-transformed matrices $\tilde{\boldsymbol{C}}$:

$$p_{\tilde{\boldsymbol{C}}}(\tilde{\boldsymbol{C}}|\tilde{\boldsymbol{\Sigma}}) = \frac{L^{LD}f(\tilde{\boldsymbol{C}})}{\Gamma_D(L)}e^{L\operatorname{tr}(\tilde{\boldsymbol{C}}-\tilde{\boldsymbol{\Sigma}})}\exp\left(-L\operatorname{tr}(e^{\tilde{\boldsymbol{C}}}e^{-\tilde{\boldsymbol{\Sigma}}})\right) \tag{29}$$

with $f(\tilde{\boldsymbol{C}}) = |\mathrm{J}_{\exp}(\tilde{\boldsymbol{C}})| / \exp[D\operatorname{tr}(\tilde{\boldsymbol{C}})]$ a normalization factor that involves the Jacobian of log transform $|\mathrm{J}_{\exp}(\tilde{\boldsymbol{C}})|$ and is equal to 1 when $D = 1$. We call such a distribution the Wishart-Fisher-Tippett distribution, denoted as $\mathcal{WFT}(\tilde{\boldsymbol{\Sigma}}; L)$, as it generalizes the Fisher-Tippett distribution to the case where $D > 1$. The expectation and variance of $\tilde{\boldsymbol{C}}$ do not seem to be known in closed form in the literature. Nevertheless, according to [2] its trace (which coincides with the logarithm determinant of $\boldsymbol{C}$: $\operatorname{tr}\tilde{\boldsymbol{C}} = \log|\boldsymbol{C}|$) has the following statistics

$$\operatorname{tr}\tilde{\boldsymbol{C}} = \operatorname{tr}\tilde{\boldsymbol{\Sigma}} + \operatorname{tr}\tilde{\boldsymbol{S}} , \tag{30}$$

$$\mathbb{E}[\operatorname{tr}\tilde{\boldsymbol{C}}] = \operatorname{tr}\tilde{\boldsymbol{\Sigma}} + \sum_{i=1}^{D}\Psi(0, L-i+1) - D\log L , \tag{31}$$

and $$\operatorname{Var}[\operatorname{tr}\tilde{\boldsymbol{C}}] = \sum_{i=1}^{D}\Psi(1, L-i+1) , \tag{32}$$

where $\tilde{\boldsymbol{S}} \sim \mathcal{WFT}(\boldsymbol{0}; L)$. This shows that the trace of the matrix logarithm suffers from additive non-zero-mean signal-independent corruptions. Note that (30) is a direct consequence of (20). If follows that the $D^2$ channels of $\tilde{\boldsymbol{C}}$ can reasonably be assumed to have approximately a stabilized variance (see Sec. III-G for numerical evidence), which opens the door to regularization with iterative schemes. Note that we cannot use the matrix log transform directly as a variance stabilization procedure, as done in the univariate homomorphic case, because we do not have an inversion formula (i.e., a bias correction formula). We will thus adopt a variational strategy.

### D. Log-channel decomposition

As mentioned in the previous paragraph, $\tilde{\boldsymbol{\Sigma}}$ lies in the vector space of complex Hermitian matrices that is isomorphic to $\mathbb{R}^{D^2}$. In this section, we describe a re-parameterization of the log-transformed covariance matrix $\tilde{\boldsymbol{\Sigma}}$ as a vector of $D^2$ reals denoted $\boldsymbol{x}$. We first define in (26) a unitary transform $\mathcal{K}$ that maps real-valued vectors $\boldsymbol{\alpha}$ of $\mathbb{R}^{D^2}$ to Hermitian $D \times D$ matrices. We also introduce a whitening affine map $(\boldsymbol{A}, \boldsymbol{b})$ so that the $D^2$ channels of $\boldsymbol{x}$ are better decorrelated, and a scaling transform $\Phi$ (i.e., a diagonal matrix) to balance the variance of noise between channels. The transform $\Omega$ between (log-transformed) covariance matrices $\tilde{\boldsymbol{\Sigma}}$ and the vector of parameter $\boldsymbol{x}$ is:

$$\tilde{\boldsymbol{\Sigma}} = \Omega(\boldsymbol{x}) = \mathcal{K}(\boldsymbol{A}\Phi\boldsymbol{x} + \boldsymbol{b}) . \tag{33}$$

---

[3]the precise definition is $\boldsymbol{\Sigma} = \sum_{n=1}^{\infty} \frac{\tilde{\boldsymbol{\Sigma}}^n}{!n}$ but is equivalent to (27) for Hermitian positive definite matrices.

[4]the matrix logarithm is not to be confused with the log-determinant function $\boldsymbol{\Sigma} \mapsto \log|\boldsymbol{\Sigma}|$, the latter equals to the trace of the former.

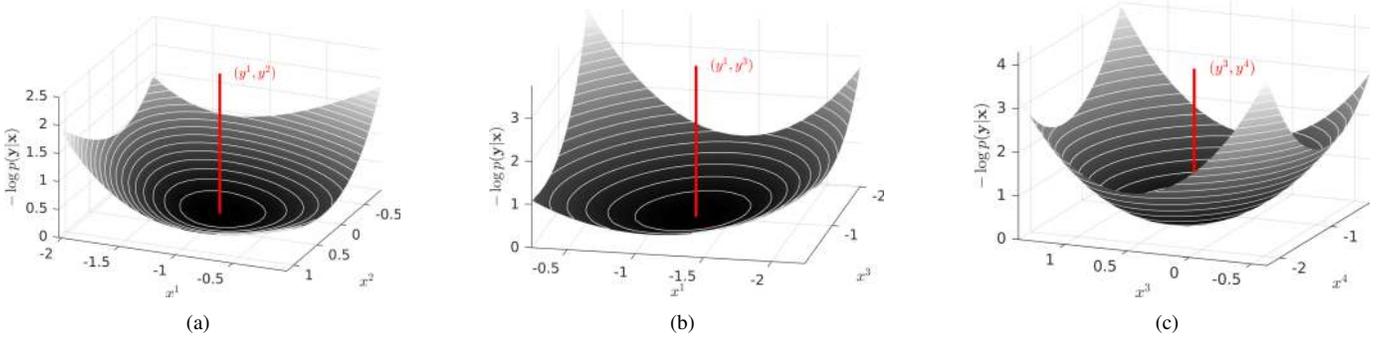

Fig. 3. Three two-dimensional sections of the Wishart-Fisher-Tippett negative log likelihood with respect to $(x^1, x^2)$, $(x^1, x^3)$ and $(x^3, x^4)$.

We denote its inverse $\Omega^{-1}(\tilde{\boldsymbol{\Sigma}}) = \Phi^{-1}\boldsymbol{A}^{-1}(\mathcal{K}^{-1}(\tilde{\boldsymbol{\Sigma}}) - \boldsymbol{b})$ and introduce the linear operator $\Omega^* = \Phi\boldsymbol{A}^*\mathcal{K}^*$, where $*$ denotes the adjoint. The re-parameterized log-transformed data are noted $\boldsymbol{y} \in \mathbb{R}^{D^2}$ and defined by $\boldsymbol{y} = \Omega^{-1}(\tilde{\boldsymbol{C}})$, see Fig. 1.

The affine transform $(\boldsymbol{A}, \boldsymbol{b})$ is chosen such that the channels $x^i$ of the re-parameterized vector are decorrelated and can then be assumed close to independent:

$$p_{\tilde{\boldsymbol{\Sigma}}}(\tilde{\boldsymbol{\Sigma}}) = p_{\boldsymbol{x}}(\boldsymbol{x}) \approx \prod_{i=1}^{D^2} p_{x^i}(x^i) . \tag{34}$$

The affine transform $(\boldsymbol{A}, \boldsymbol{b})$ is obtained by principal component analysis, as described in Sec. III-G. Moreover, the scaling matrix $\Phi$ is chosen such that the variance of $y^i$ is equal to 1 for all channels $i$, see Sec. III-G.

### E. Proposed variational approach

We now extend $\boldsymbol{x}$ and $\boldsymbol{y}$ into images of $n$ real-valued vectors each of dimension $D^2$. Writing $-\log p_{x^i}(x^i) = \mathcal{R}(x^i)$ and using the relations (34) within the MAP framework leads to the following minimization problem

$$\hat{\boldsymbol{x}} \in \underset{\boldsymbol{x}}{\operatorname{argmin}} \; -\log p_{\boldsymbol{y}}(\boldsymbol{y}|\boldsymbol{x}) + \sum_{i=1}^{D^2} \mathcal{R}(x^i) , \tag{37}$$

where, from (29), we have

$$-\log p_{\boldsymbol{y}}(\boldsymbol{y}|\boldsymbol{x}) = L\sum_{k=1}^{n} \operatorname{tr}\left(\Omega(\boldsymbol{x}_k) + e^{\Omega(\boldsymbol{y}_k)}e^{-\Omega(\boldsymbol{x}_k)}\right) + \text{Cst.} , \tag{38}$$

and the final estimate $\hat{\boldsymbol{\Sigma}}_k$ at pixel index $k$ is defined as $\hat{\boldsymbol{\Sigma}}_k = \exp\Omega(\hat{\boldsymbol{x}}_k)$. Remark that this problem boils down to (10) when $D = 1$. Note also the difference in notations between $\boldsymbol{x}_k$ a vector of $D^2$ coefficients at pixel $k$ and $x^i$ a scalar image of $n$ pixels corresponding to the $i$-th channel of $\boldsymbol{x}$.

While the direct multivariate variational approach (25) requires optimizing on the cone of complex Hermitian matrices, a first major advantage is that Problem (37) is unconstrained on $\mathbb{R}^{D^2}$. Even though the likelihood term (38) is still nonconvex, it appears to be much more suitable for optimization with less local minima. Figure 3 illustrates the evolution of this term (with $D = 2$, and $n = 1$) on three two-dimensional cross sections of $\mathbb{R}^4$ showing that for each of such sections the energy appears to be convex. Nevertheless, note that convexity in low-dimensional spaces does not necessarily guarantee convexity in higher dimensional spaces. Unlike the direct multivariate variational approach (25), (38) appears to be convex in many scenarios, e.g., in the univariate case (10), and more generally when $\hat{\boldsymbol{x}}_k$ is restricted to the vectorial subspace where $\Omega(\hat{\boldsymbol{x}}_k)$ commutes with $\Omega(\boldsymbol{y}_k)$, hence satisfying $e^{\Omega(\boldsymbol{y}_k)}e^{-\Omega(\hat{\boldsymbol{x}}_k)} = e^{\Omega(\boldsymbol{y}_k - \hat{\boldsymbol{x}}_k)}$. Convexity, in this case, follows from the fact that the trace is linear, hence convex, the function $\hat{\boldsymbol{x}}_k \to \operatorname{tr} e^{\Omega(-\hat{\boldsymbol{x}}_k + \boldsymbol{y}_k)}$ is a convex spectral function [27], [43], and $\Omega$ is affine.

As for the mono-dimensional case, we will thus consider the ADMM algorithm which iteratively performs the updates

$$\hat{\boldsymbol{z}} \leftarrow \underset{\boldsymbol{z}}{\operatorname{argmin}} \; \frac{\beta}{2}\|\boldsymbol{z} - \hat{\boldsymbol{x}} + \hat{\boldsymbol{d}}\|^2 + \sum_{i=1}^{D^2} \mathcal{R}(z^i) , \tag{39}$$

$$\hat{\boldsymbol{d}} \leftarrow \hat{\boldsymbol{d}} + \hat{\boldsymbol{z}} - \hat{\boldsymbol{x}} , \tag{40}$$

$$\hat{\boldsymbol{x}} \leftarrow \underset{\boldsymbol{x}}{\operatorname{argmin}} \; \frac{\beta}{2}\|\boldsymbol{x} - \hat{\boldsymbol{z}} - \hat{\boldsymbol{d}}\|^2 - \log p_{\boldsymbol{y}}(\boldsymbol{y}|\boldsymbol{x}) . \tag{41}$$

Since the noise variance is approximately stabilized and close to 1 on each channel (see, III-G), a single value for $\beta$ can be chosen for all channels. With such a choice, the same amount of smoothing is reached in practice in all regions after a finite number of iterations.

Equation (39) is separable on the different channels, and can then be solved by applying $D^2$ times a regularized least-square solver. Equation (41) amounts to solving $n$ separable problems of the form

$$\underset{\boldsymbol{x}_k}{\operatorname{argmin}} \; \frac{\beta}{2}\|\boldsymbol{x}_k - \boldsymbol{a}_k\|^2 + L\operatorname{tr}\left(\Omega(\boldsymbol{x}_k) + e^{\Omega(\boldsymbol{y}_k)}e^{-\Omega(\boldsymbol{x}_k)}\right) \tag{42}$$

where $\boldsymbol{a}_k = \hat{\boldsymbol{z}}_k + \hat{\boldsymbol{d}}_k$. As for the univariate case, we will consider Newton's method to solve this optimization problem. Its gradient is given by

$$\beta(\boldsymbol{x}_k - \boldsymbol{a}_k) +$$
$$L\Omega^*\left(\mathbf{Id} - \int_0^1 e^{(u-1)\Omega(\boldsymbol{x}_k)}e^{\Omega(\boldsymbol{y}_k)}e^{-u\Omega(\boldsymbol{x}_k)}\mathrm{d}u\right) , \tag{43}$$

$$\Delta_\infty^{k,i} = \frac{\beta(x_k^i - a_k^i) + L\left[\Omega^*\left(\mathbf{Id} - \int_0^1 e^{(u-1)\Omega(\boldsymbol{x}_k)}e^{\Omega(\boldsymbol{y}_k)}e^{-u\Omega(\boldsymbol{x}_k)}\mathrm{d}u\right)\right]^i}{\left|\beta + L\Omega^*\int_0^1 e^{(u-1)\Omega(\boldsymbol{x}_k)}e^{\Omega(\boldsymbol{y}_k)}e^{-u\Omega(\boldsymbol{x}_k)}\mathrm{d}u\right|^i} . \tag{35}$$

$$\Delta_Q^{k,i} = \frac{\beta(x_k^i - a_k^i) + L\left[\Omega^*\left(\mathbf{Id} - \frac{1}{Q}\sum_{q=1}^Q e^{(u_q-1)\Omega(\boldsymbol{x}_k)}e^{\Omega(\boldsymbol{y}_k)}e^{-u_q\Omega(\boldsymbol{x}_k)}\right)\right]^i}{\left|\beta + L\Omega^*\frac{1}{Q}\sum_{q=1}^Q e^{(u_q-1)\Omega(\boldsymbol{x}_k)}e^{\Omega(\boldsymbol{y}_k)}e^{-u_q\Omega(\boldsymbol{x}_k)}\right|^i} , \quad u_q = \frac{q - \frac{1}{2}}{Q} . \tag{36}$$

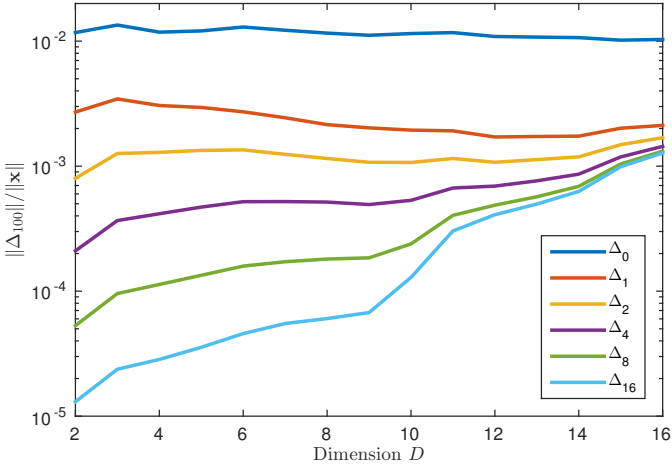

Fig. 4. Evolution with respect to the matrix dimension $D$ of the error $\|\Delta_{100}\|/\|\hat{\boldsymbol{x}}\|$ after 10 iterations of Newton's method for which the update $\Delta_\infty$ is replaced by $\Delta_Q$, for $Q = \{0, 1, 2, 4, 8, 16\}$. For each $D$, the results are averaged over 100 realizations of $\boldsymbol{C}$ (with $L = D$), each realization is obtained from a different version of $\boldsymbol{\Sigma}$ (also randomly generated).

see the development in Appendix A. By mimicking the univariate case, we consider the following approximation for the Hessian

$$\mathrm{diag}\left(\beta\mathbf{Id} + L\Omega^*\int_0^1 e^{(u-1)\Omega(\boldsymbol{x}_k)}e^{\Omega(\boldsymbol{y}_k)}e^{-u\Omega(\boldsymbol{x}_k)}\mathrm{d}u\right), \tag{44}$$

leading to the quasi-Newton iteration

$$\hat{x}_k^i \leftarrow \hat{x}_k^i - \Delta_\infty^{k,i} , \tag{45}$$

where $\Delta_\infty^{k,i}$ is defined in eq. (35). Note that the expression becomes the exact Newton iteration when $\Omega(\boldsymbol{x}_k)$ and $\Omega(\boldsymbol{y}_k)$ are commuting matrices. As for the univariate case, we notice that ten iterations is enough to offer good performance with a reasonable computational time.

The update matrix $\Delta_Q^{k,i}$, defined in eq. (35), requires numerical integration. We use Riemann integral approximation with $Q$ rectangles as defined in eq. (36). Of course, the smaller $Q$, the faster the iteration (45). Figure 4 shows the evolution of $\|\Delta_{100}\|/\|\hat{\boldsymbol{x}}\|$ with respect to $D$ (with $L = D$ and $\beta = 10L$) after 10 iterations of Newton's method where $\Delta_\infty$ is substituted by $\Delta_Q$, for $Q = \{1, 2, 4, 8, 16\}$. We furthermore define $\Delta_0$ by substituting the integral by $\exp(\Omega(\boldsymbol{y}_k - \boldsymbol{x}_k))$ as a crude approximation of the gradient (and exact in the case where $\Omega(\hat{\boldsymbol{x}}_k)$ commutes with $\Omega(\boldsymbol{y}_k)$). Ideally, the algorithm would reach the optimum solution such that $\Delta_\infty \approx \Delta_{100} = 0$ (as for the univariate case). This error is quite small whatever $Q$ and $D$. For low values of $D$, the error is reduced by about a factor 10 each time $Q$ is multiplied by 4. The increase of the error with $D$ is due to an accumulation of numerical errors of the exponential matrix function, all the more important when $Q$ gets larger. This experiment reveals that $Q = 1$ is a good choice in practice, reaching a relative error of about $2 \cdot 10^{-3}$ whatever $D \in [2, 16]$ while requiring a reasonable computation time.

### F. Adaptation to advanced filters

As in the univariate case, we notice that ADMM converges when the regularizer $\mathcal{R}$ is chosen as nonconvex. More remarkably, as observed in many contexts, *e.g.*, [67], [13], [65], [14], [56], [8], replacing the minimization problem in (39) by the solution of a Gaussian denoiser – a scheme known as plug-and-play – leads to appealing results. The resulting algorithm, coined MUlti-channel LOgarithm with Gaussian denoising (MuLoG), is given as

$$\hat{z}^i \leftarrow f_{\beta^{-1/2}}(\hat{x}^i - \hat{d}^i), \quad \text{for } i = 1, \ldots, D^2 , \tag{46}$$

$$\hat{\boldsymbol{d}} \leftarrow \hat{\boldsymbol{d}} + \hat{\boldsymbol{z}} - \hat{\boldsymbol{x}} , \tag{47}$$

$$\hat{\boldsymbol{x}} \leftarrow \underset{\boldsymbol{x}}{\mathrm{argmin}} \frac{\beta}{2}\|\boldsymbol{x} - \hat{\boldsymbol{z}} - \hat{\boldsymbol{d}}\|^2 - \log p_{\boldsymbol{y}}(\boldsymbol{y}|\boldsymbol{x}) , \tag{48}$$

where $f_{\sigma^2}$ is again a denoiser for images contaminated by zero-mean white Gaussian noise $\mathcal{N}(0; \sigma^2\mathrm{Id})$, and (48) is solved as mentioned in the previous paragraph. This plug-and-play ADMM algorithm[5] is proven to converge[6] [9] as soon as $f_{\sigma^2}$ is a *bounded* denoiser, i.e., satisfying for all $x$

$$\|f_{\sigma^2}(x) - x\|^2 \leq n\sigma^2 C \tag{49}$$

for some constant $C$ independent of $n$ and $\sigma$. As mentioned in [9], this is a weak condition which can be expected from most denoisers. In particular, this condition implies that $f_{\sigma^2}$ should tend to the identity when $\sigma$ tends to 0. We depict this algorithm on the last row of Fig.1.

---

[5] this requires to add a small update of $\beta$ during the iterations, according to the rule given in [9].

[6] convergence holds in our case by continuity of the gradient of (37), which implies that the gradient is Lipschitz in any compact subsets of $\mathbb{R}^{D^2 \times n}$.

## G. Calibration

In practice, the operators $\boldsymbol{A}$, $\boldsymbol{b}$ and $\Phi$ are obtained from the log matrix of the input image $\boldsymbol{C}$ with $n$ pixels. Let $\{\boldsymbol{\alpha}_1, \ldots, \boldsymbol{\alpha}_n\}$ be the collection of $n$ vectors extracted from the matrix logarithm of the input image as $\boldsymbol{\alpha}_k = \mathcal{K}^{-1}(\tilde{\boldsymbol{C}}_k)$. We thus define $\boldsymbol{b} = \frac{1}{n} \sum_{k=1}^{n} \boldsymbol{\alpha}_k$ and $\boldsymbol{A}$ the matrix whose columns are the eigenvectors (with unit norm) of $\frac{1}{n} \sum_{k=1}^{n} (\boldsymbol{\alpha}_k - \boldsymbol{b})(\boldsymbol{\alpha}_k - \boldsymbol{b})^t$. As a result, the vector $\boldsymbol{A}^{-1}(\boldsymbol{\alpha}_k - \boldsymbol{b})$ is the representation of $\boldsymbol{\alpha}_k$ in the principal component analysis of $\{\boldsymbol{\alpha}_1, \ldots, \boldsymbol{\alpha}_n\}$ known to maximize inter-channel decorrelation, and thus leads to (34). By assuming the noise variance is stabilized on each channel of $\boldsymbol{A}^{-1}(\boldsymbol{\alpha} - \boldsymbol{b})$, we define the diagonal matrix $\Phi$ as $\mathrm{diag}(\hat{\sigma}_1^2, \ldots, \hat{\sigma}_D^2)$, where $\hat{\sigma}_1^2, \ldots, \hat{\sigma}_D^2$ are obtained by estimating separately the variance of the noise on each channel of $\boldsymbol{A}^{-1}(\boldsymbol{\alpha} - \boldsymbol{b})$. This procedure ensures that the noise variance is about 1 on each channel of $\boldsymbol{y}$. In practice, we resort to the median absolute deviation estimate (MAD) [33] for this task, but other strategies could be used, such as [45], [46], [62].

Figure 5 displays the nine components of $\boldsymbol{y}$ for an image of $3 \times 3$ matrices $\boldsymbol{C}$. While the noise variance in $\boldsymbol{C}$ is clearly signal-dependent, the amplitudes of the fluctuations in the different channels of $\boldsymbol{y}$ appear fairly constant whatever the underlying signal $\boldsymbol{\Sigma}$. Nevertheless, the noise is not signal-independent, we can observe that some regions are a smidgen noisier than others. Remark also that the noise variance is also about the same for all channels. Then, we can reasonably claim from this experiment that our log channel decomposition approximately stabilizes the noise variance.

Thanks to this approximate stabilization, we notice that choosing 6 iterations with $\beta = 1 + 2/L$ provides satisfying solutions in all tested situations, whatever the dimension $D$, the number of looks $L$ and the embedded Gaussian denoiser. The initialization is chosen as $\hat{\boldsymbol{x}} = \boldsymbol{y}$, $\hat{z}^i = f_1(y^i)$ and $\hat{\boldsymbol{d}} = \hat{\boldsymbol{z}} - \hat{\boldsymbol{x}}$. We always choose $Q = 1$ according to Sec. III-E.

When $L < D$, the matrix $\boldsymbol{C}_k$ is rank deficient, and the likelihood term (38), though convex when restricted to commutative matrices, is however no longer strictly convex and has an infinite number of minimizers. Worse, its logarithm $\tilde{\boldsymbol{C}}_k$ is undefined which prevents the computation of $\boldsymbol{A}$, $\boldsymbol{b}$ and $\Phi$ used in our log-channel decomposition. To deal with these issues a practical solution is to perform a small diagonal loading of the input matrices $\boldsymbol{C}_k$ to enforce their positive definite property. After trying several other alternatives, we found that satisfying results were reached by performing a spatially varying re-scaling of the off-diagonal elements as

$$(\boldsymbol{C}_k)_{i,j}^{\text{new}} \leftarrow \frac{|\sum_l w_{k,l}(\boldsymbol{C}_l)_{i,j}|}{\sqrt{\sum_l w_{k,l}(\boldsymbol{C}_l)_{ii} \sum_l w_{k,l}(\boldsymbol{C}_l)_{jj}}} (\boldsymbol{C}_k)_{i,j} , \quad (50)$$

where $w_{k,l}$ are the weights of a Gaussian kernel with bandwidth 1 pixel. This procedure inevitably introduces some bias in the solution, but this bias does not seem to be significant.

## IV. Efficient matrix log and exp transforms

Our algorithm requires to compute several times $n$ matrix logarithms, matrix exponentials, and matrix-by-matrix products. It is thus important to make these operations as fast as

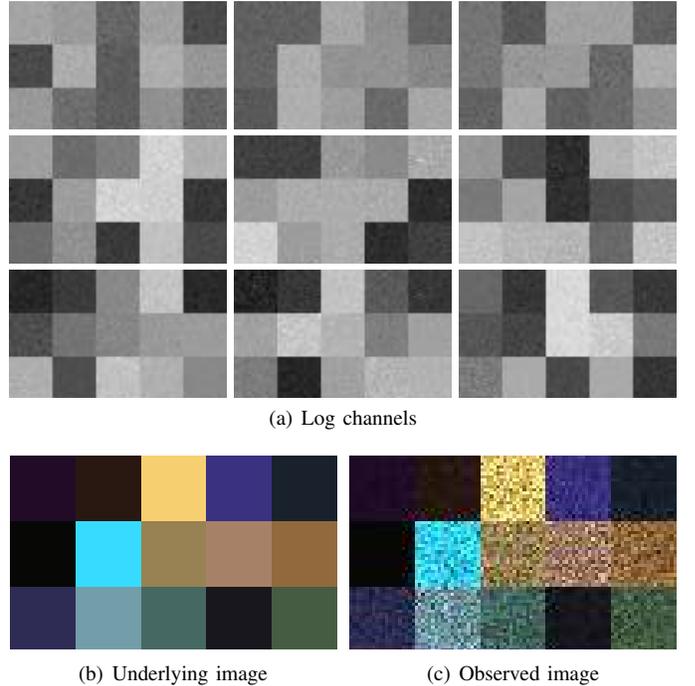

(a) Log channels

(b) Underlying image      (c) Observed image

Fig. 5. (a) the nine log channels $\boldsymbol{x}^1, \ldots, \boldsymbol{x}^9$, with approximate unit stabilized variance, of (c) the image $\boldsymbol{C}$ generated from (b) the image $\boldsymbol{\Sigma}$ displayed with a RGB representation based on its decomposition in Pauli's basis.

possible. Our first attempt to call $n$ times Matlab's functions `logm`, `expm` and `mtimes` leads to very slow computations. It is, in fact, important to use a vectorial implementation of these functions. When $D = 2$, we can write the matrix $\boldsymbol{C}$ and $\tilde{\boldsymbol{C}}$ as

$$\boldsymbol{C} = \begin{pmatrix} a & c^* \\ c & b \end{pmatrix} \quad \text{and} \quad \tilde{\boldsymbol{C}} = \begin{pmatrix} \tilde{a} & \tilde{c}^* \\ \tilde{c} & \tilde{b} \end{pmatrix} , \quad (51)$$

and the relation between both are given, for $\tilde{\boldsymbol{C}} = \log \boldsymbol{C}$, as

$$\begin{cases} \tilde{a} &= \left[(a - b + \delta)\ell_1 - (a - b - \delta)\ell_2\right]/(2\delta) , \\ \tilde{b} &= \left[(b - a + \delta)\ell_1 - (b - a - \delta)\ell_2\right]/(2\delta) , \\ \tilde{c} &= c(\ell_1 - \ell_2)/\delta , \end{cases} \quad (52)$$

and, for $\boldsymbol{C} = \exp \tilde{\boldsymbol{C}}$, as

$$\begin{cases} a &= \left[(\tilde{a} - \tilde{b} + \tilde{\delta})e_1 - (\tilde{a} - \tilde{b} - \tilde{\delta})e_2\right]/(2\tilde{\delta}) , \\ b &= \left[-(\tilde{a} - \tilde{b} - \tilde{\delta})e_1 - (\tilde{a} - \tilde{b} + \tilde{\delta})e_2\right]/(2\tilde{\delta}) , \\ c &= \tilde{c}(e_1 - e_2)/\tilde{\delta} , \end{cases} \quad (53)$$

where $\delta = \sqrt{4|c|^2 + (a - b)^2}$, $\ell_1 = \log\left[(a + b + \delta)/2\right]$, $\ell_2 = \log\left[(a + b - \delta)/2\right]$, $\tilde{\delta} = \sqrt{4|\tilde{c}|^2 + (\tilde{a} - \tilde{b})^2}$, $e_1 = \exp\left[(\tilde{a} + \tilde{b} + \tilde{\delta})/2\right]$ and $e_2 = \exp\left[(\tilde{a} + \tilde{b} - \tilde{\delta})/2\right]$. With these established relations, we can compute the matrix logarithm and matrix exponential in a vectorial way, starting by computing for the $n$ pixels the coefficient $\tilde{a}$, and next $\tilde{b}$ and finally $\tilde{c}$. We have implemented this vectorial procedure in Matlab (using the element-wise operators `.*` and `./`). and observed an acceleration by a factor larger than 600 compared to calling $n$ times the dedicated Matlab's function. We were also able to derive closed form formula for the case where $D = 3$,

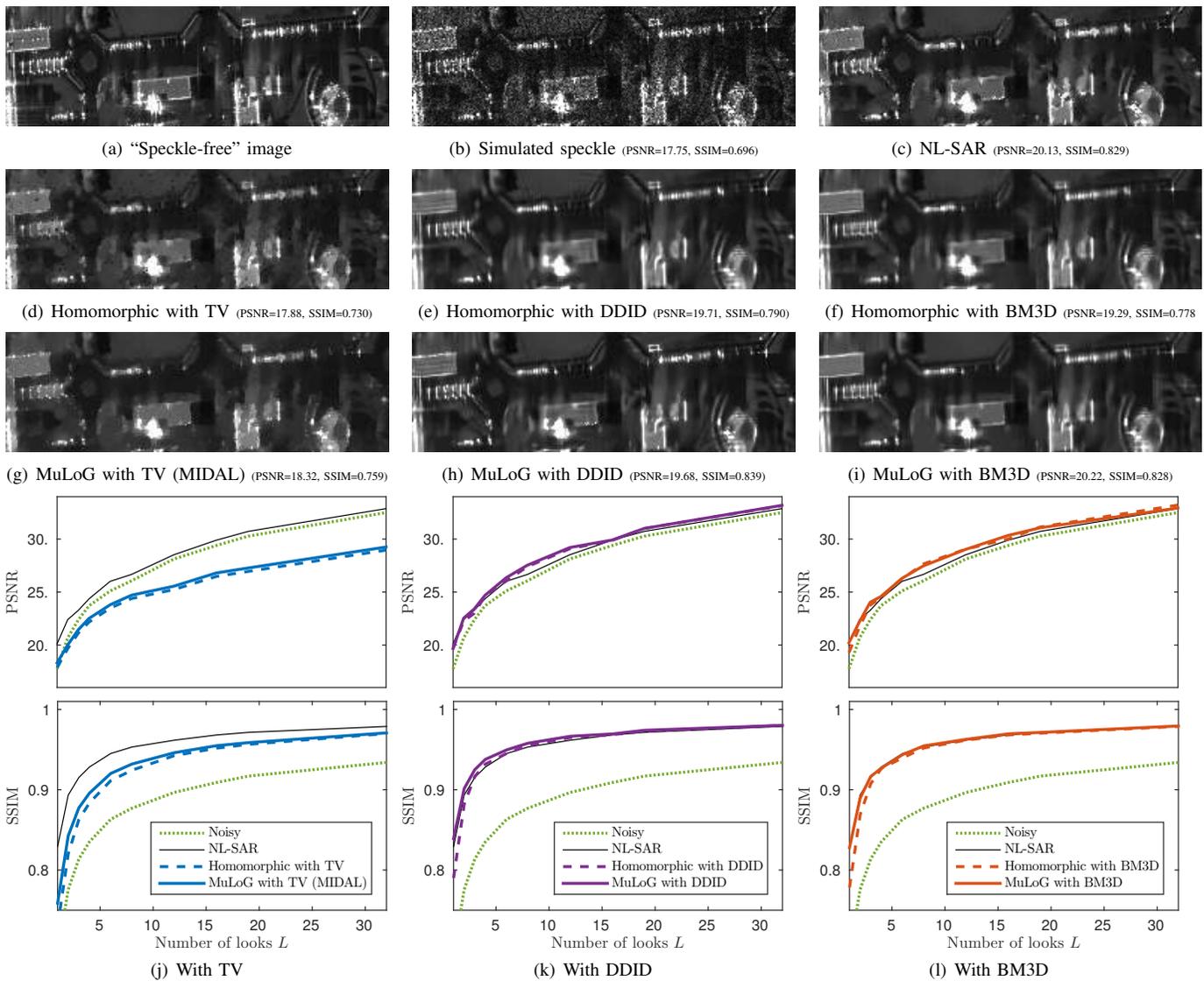

(a) "Speckle-free" image    (b) Simulated speckle (PSNR=17.75, SSIM=0.696)    (c) NL-SAR (PSNR=20.13, SSIM=0.829)

(d) Homomorphic with TV (PSNR=17.88, SSIM=0.730)    (e) Homomorphic with DDID (PSNR=19.71, SSIM=0.790)    (f) Homomorphic with BM3D (PSNR=19.29, SSIM=0.778)

(g) MuLoG with TV (MIDAL) (PSNR=18.32, SSIM=0.759)    (h) MuLoG with DDID (PSNR=19.68, SSIM=0.839)    (i) MuLoG with BM3D (PSNR=20.22, SSIM=0.828)

(j) With TV    (k) With DDID    (l) With BM3D

Fig. 6. (a) A "speckle-free" SAR image, (b) its simulated version with single-look speckle ($L = 1$), Results of speckle reduction with (c) NL-SAR [20], (d-f) MuLoG embedding respectively TV, DDID and BM3D, (g-i) the homomorphic approach embedding respectively TV, DDID and BM3D. (j-l) Measure of performance on this image in terms of PSNR (the larger the better) and SSIM (the larger the better) metrics as a function of the number of looks $L$.

leading in this case to an acceleration by a factor 75. The formulas being too long to be inserted in this paper, we invite the interested reader to refer to [16]. The same workaround is applied for the matrix-by-matrix product.

## V. NUMERICAL EXPERIMENTS

*a) Considered Gaussian denoisers:* In this section, numerical experiments are conducted based on three Gaussian denoisers. First, we consider the classical isotropic total-variation (TV) [58], i.e., $\mathcal{R}(x^i) = \lambda \sum_{k=1}^{n} \|(\nabla x^i)_k\|$ which corresponds to the multivariate extension of [7]. The parameter $\lambda$ has been tuned to 0.7 and kept the same for all data-sets, whatever the dimension $D$ or the number of looks $L$. This can be achieved only because noise components in the different channels have all been stabilized to a variance of 1. Next, we consider DDID [38], a recent hybrid method based on bilateral filtering and shrinkage of Fourier coefficients. Finally, we consider BM3D [12], a ten-years algorithm, reaching remarkable results with fast computation, based on patch similarity and three-dimensional wavelet shrinkage. The inner parameters of these two latter algorithms have been kept the same as the one provided in the authors' implementation.

*b) Simulations in the univariate case:* Figure 6 reports the results of a simulated experiment in the SAR amplitude context ($D = 1$). Speckle has been generated on top of a $409 \times 409$ "speckle-free" SAR image[7] of the CNES in Toulouse, France (obtained by ONERA's RAMSES sensor). These images having large dynamic ranges, they have been saturated and gamma-corrected for display purposes. Speckle with levels from $L = 1$ to 32 have been considered. The performance of MuLoG is compared to the homomorphic approach when both embed the aforementioned Gaussian de-

---

[7] Speckle in this 1 meter ground-resolution image has been reduced by a factor 100 by averaging a 10 centimeters ground-resolution image.

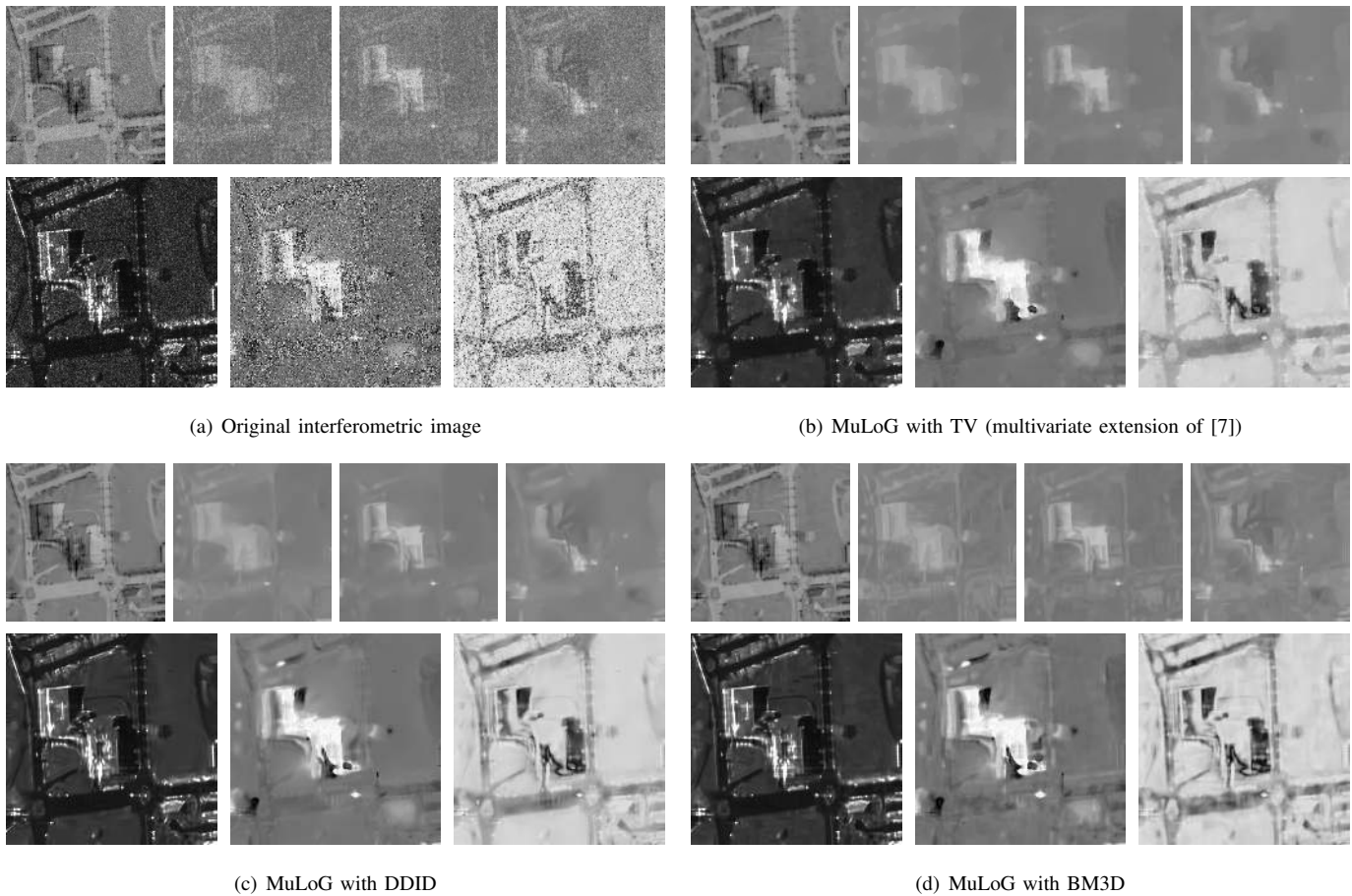

**(a)** Original interferometric image

**(b)** MuLoG with TV (multivariate extension of [7])

**(c)** MuLoG with DDID

**(d)** MuLoG with BM3D

Fig. 7. (a) A single-look interferometric image ($D = 2$, $L = 1$) of Toulouse (France) sensed by RAMSES (©ONERA). Estimation with our proposed multivariate framework MuLoG for (b) Total-Variation (TV), (c) DDID and (d) BM3D. For each image, the four log channels are displayed on top $\boldsymbol{x}^1, \boldsymbol{x}^2, \boldsymbol{x}^3, \boldsymbol{x}^4$, below are the trace $\operatorname{tr} \boldsymbol{C}$, the phase $\arg \boldsymbol{C}_{12}$ and the coherence $|\boldsymbol{C}_{12}|/\sqrt{|\boldsymbol{C}_{11}\boldsymbol{C}_{22}|}$ respectively.

noisers. As a baseline, the performance are also compared to a recent patch-based speckle reduction technique: NL-SAR [20]. In the single-look case ($L = 1$), close-ups of size $256 \times 80$ are given for visual inspection. PSNR is first used as a quantitative measure of amplitude reconstruction quality[8]. In this high dynamic range setting, PSNR values are extremely sensitive to a small discrepancy, all the more for tiny bright structures. As a consequence, it may appear irrelevant, for instance, by ranking MIDAL worst than the noisy image itself[9]. For this reason, we have also included SSIM [66] values. Visual inspection and quantitative criteria seem to indicate that MuLoG provides more relevant solutions than the homomorphic approach. In particular, the homomorphic approach tends to oversmooth bright targets and leaves residual dark stains, revealing its lack of robustness against the heavy left tail of the Fisher-Tippet distribution. Overall, results obtained with MuLoG are on a par with those of NL-SAR, while displaying small variations depending on the chosen Gaussian denoiser. Visual analysis of the output produced by each variant of MuLoG may be useful to discard structures appearing only with one Gaussian denoiser as artifacts.

*c) Results in the multivariate case:* Figures 7, 8 and 9 give three illustrations of MuLoG in real speckle reduction contexts. Figure 7 corresponds to a $256 \times 256$ airborne single-look SAR interferometric image ($D = 2$, $L = 1$) of a building in Toulouse, France (sensed by RAMSES). Figure 8 corresponds to a $250 \times 250$ spaceborne single-look SAR interferometric image ($D = 2$, $L = 1$) of the dam of Serre-Ponçon, France (sensed by TerraSAR-X). Figure 9 corresponds to a $512 \times 512$ airborne SAR polarimetric image ($D = 3$, $L = 1$) of the city of Kaufbeuren, Germany (sensed by F-SAR).

In the sub-figures (b) are displayed the results for the isotropic total-variation (TV) regularization. In this case, ADMM finds a local minimum of (37). We can observe that the solution in each channel inherits from the well-known behavior of the univariate TV regularization: the solution appears piece-wise constant, small details and low contrasted features are lost, and a slight bias can be observed.

In the sub-figures (c) and (d) are displayed the results for DDID [38] and BM3D [12]. In these cases, our proposed algorithm does not explicitly minimize a given energy but we can observe that it converges in practice to relevant solutions. As for TV, we can observe that the solution in each channel inherits from the behavior of the original method: in this case

---

[8] Since images are not 8bit, the standard PSNR has been modified with a peak value defined as the 99th quantile of the "speckle-free" SAR image.

[9] TV regularization produces a systematic loss of contrast (see, [61], [21]), even stronger for punctual bright structures.

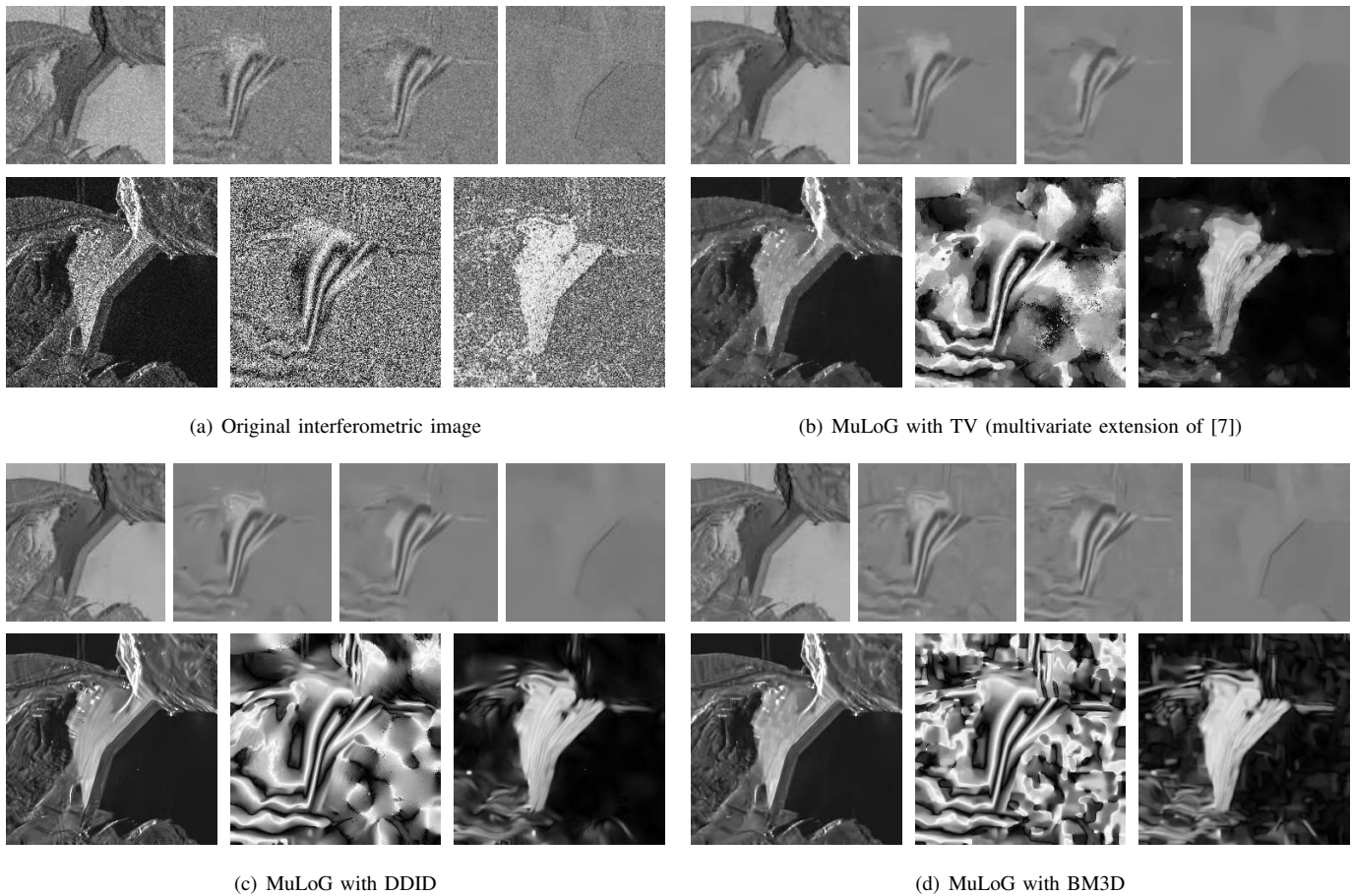

(a) Original interferometric image

(b) MuLoG with TV (multivariate extension of [7])

(c) MuLoG with DDID

(d) MuLoG with BM3D

Fig. 8. (a) A single-look interferometric image ($D = 2, L = 1$) of Serre-Ponçon (France) sensed by TerraSar-X (image courtesy to Airbus Defence and Space). Estimation with our proposed multi-variate framework MuLoG for (b) Total-Variation (TV), (c) DDID and (d) BM3D. For each image, the four log channels are displayed on top $x^1, x^2, x^3, x^4$, below are the trace $\mathrm{tr}\,C$, the phase $\arg C_{12}$ and the coherence $|C_{12}|/\sqrt{|C_{11}C_{22}|}$ respectively.

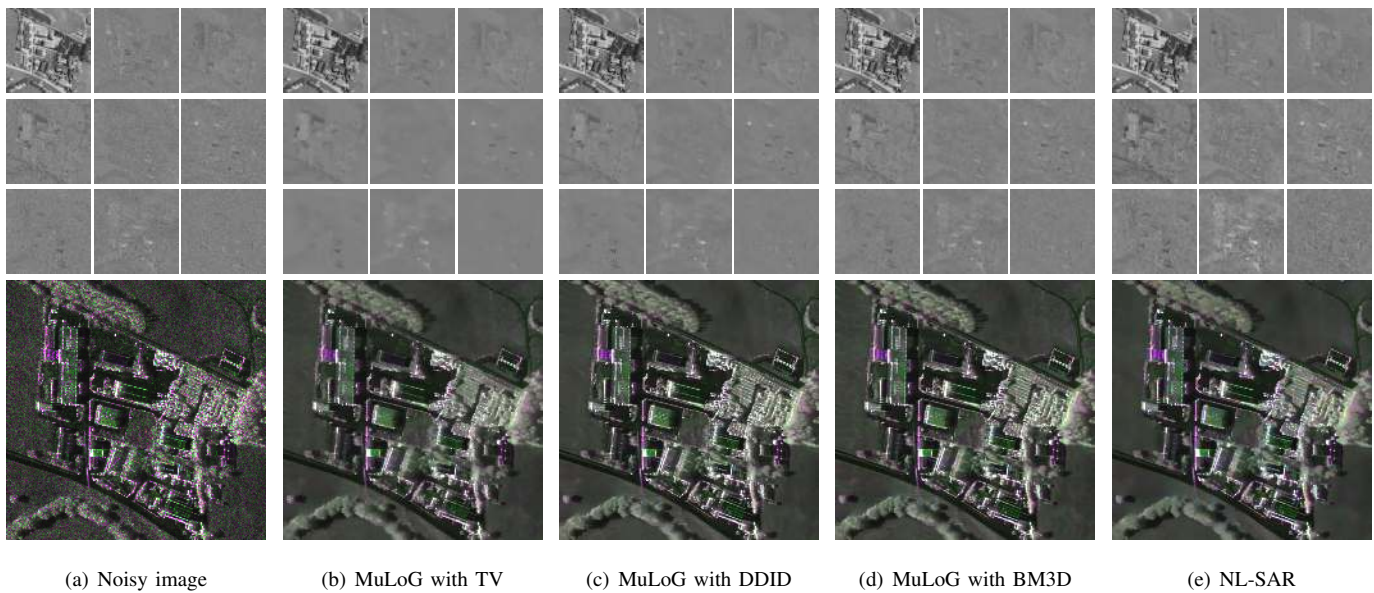

(a) Noisy image   (b) MuLoG with TV   (c) MuLoG with DDID   (d) MuLoG with BM3D   (e) NL-SAR

Fig. 9. (a) A single-look interferometric image ($D = 3, L = 1$) of Kaufbeuren (Germany) sensed by F-SAR (©DLR). Estimation with (b) NL-SAR [20] our proposed multi-variate framework MuLoG for (c) Total-Variation (TV), (d) DDID and (e) BM3D. For each image, the nine log channels are displayed on top $x^1, \ldots, x^9$, below are a RGB representation of $\hat{\Sigma}$ based on its decomposition in Pauli's basis.

the small details and low contrasted features are well restored but some known typical artifacts (small oscillating features) of these two methods can be observed as well.

Finally, we display in sub-figure 9(e) the result of NL-SAR [20] in the SAR polarimetric context. It provides a baseline for assessing the relative performance of MuLoG against a recent patch-based speckle reduction technique dedicated to SAR imagery. Compared to the different MuLoG versions, NL-SAR provides slightly over-smoothed results but is free of systematical oscillating artifacts inherent of DDID and BM3D.

*d) Computation time analysis:* Table I reports the computation time of our proposed approach with different embedded Gaussian denoisers, different image sizes and different covariance matrix sizes. Our implementation uses parallelization obtained by running eq. (46) in parallel on the $D^2$ channels, and running eq. (48) on different subsets of the $n$ pixels. Our experiments were conducted on a processor with 4 cores. In this experiment, we observe that the total computation time is always faster than running successively $6 \times D^2$ times the original algorithm.

This paper only deals with the problem of applying Gaussian filters to multi-variate SAR data. A deeper study comparing their relative performance or comparing them to state-of-the-art multi-variate SAR filters is out of the scope of this paper. However, we believe that this study should be performed in a future work and is of main interest for the community. For this reason, we have released a Matlab script[10], under CECILL license, that takes as input a multivariate SAR image, its number of looks and a Gaussian denoiser provided by the user, and outputs the filtered SAR image.

## VI. Conclusion

While a large diversity of speckle reduction methods exist for intensity images, only few can be extended to process multi-channel SAR images. In particular, the extension of variational methods leads to challenging optimization problems due to the nonconvexity of the neg-log-likelihood and the positive definiteness constraint on the covariance matrices. Furthermore, signal and channel dependent variance lead to restoration results with an uneven suppression of speckle. This paper introduced a general scheme based on a matrix logarithm transform to approximately stabilize speckle variance and produce close to independent channels. Each channel can then be processed with a user-defined Gaussian denoiser. Upon re-iterating, a good fit of the restored multi-channel image with the Wishart distribution of input covariance matrices is enforced.

Special care is paid to ensure that the method requires no parameter tuning other than possibly within the Gaussian denoiser, and that these parameters, if any, can be tuned once for all. The resulting generic method can then include Gaussian denoisers selected by the user and tremendously extends the set of available speckle reduction methods applicable to multi-channel SAR images. We believe that this offers several notable advantages: (i) the SAR imaging community will directly benefit from upcoming progress made in the field of image denoising; (ii) several images with reduced speckle can be produced by very different denoising algorithms, and these images can be compared to discard artifacts and confirm weak structures; (iii) this family of speckle reduction methods can serve as a reference to benchmark future specialized speckle reduction algorithms. This motivated the release of an open-source code implementing our method[10].

## Appendix A
### Gradient of the neg log likelihood

In this section, we establish the part of formula (43) corresponding to the gradient of $\operatorname{tr}\left(\Omega(\boldsymbol{x}) + e^{\Omega(\boldsymbol{y})}e^{-\Omega(\boldsymbol{x})}\right)$. To this end, recall that for a real valued differentiable function $f$,

$$\mathrm{d}f(\boldsymbol{x}) = \operatorname{tr}[g(\boldsymbol{x})\mathrm{d}\boldsymbol{x}] \Leftrightarrow \nabla_{\boldsymbol{x}}f(\boldsymbol{x}) = g(\boldsymbol{x})^* \ . \quad (54)$$

where $g$ is a matrix valued function whose dimension depends on that of $\boldsymbol{x}$, and $*$ denotes the adjoint operator. Since $\Omega(\boldsymbol{x}) = \mathcal{K}(\boldsymbol{A}\Phi\boldsymbol{x} + \boldsymbol{b})$ is affine, this directly implies that

$$\nabla_{\boldsymbol{x}}\operatorname{tr}\Omega(\boldsymbol{x}) = \Phi^*\boldsymbol{A}^*\mathcal{K}^* = \Omega^* \ . \quad (55)$$

We now use that for any matrix valued function $h$, we have (proven in [37], [5] according to [51])

$$\mathrm{d}e^{h(\boldsymbol{x})} = \int_0^1 e^{uh(\boldsymbol{x})}(\mathrm{d}h(\boldsymbol{x}))e^{(1-u)h(\boldsymbol{x})}\mathrm{d}u \ . \quad (56)$$

It follows that, for any matrix $\boldsymbol{L}$, we have

$$\mathrm{d}\operatorname{tr}\left[\boldsymbol{L}e^{h(\boldsymbol{x})}\right] = \operatorname{tr}\left[\boldsymbol{L}(\mathrm{d}e^{h(\boldsymbol{x})})\right] \ , \quad (57)$$

$$= \operatorname{tr}\left[\boldsymbol{L}\int_0^1 e^{uh(\boldsymbol{x})}(\mathrm{d}h(\boldsymbol{x}))e^{(1-u)h(\boldsymbol{x})}\mathrm{d}u\right] \ , \quad (58)$$

$$= \int_0^1 \operatorname{tr}\left[\boldsymbol{L}e^{uh(\boldsymbol{x})}(\mathrm{d}h(\boldsymbol{x}))e^{(1-u)h(\boldsymbol{x})}\right]\mathrm{d}u \ , \quad (59)$$

$$= \int_0^1 \operatorname{tr}\left[e^{(1-u)h(\boldsymbol{x})}\boldsymbol{L}e^{uh(\boldsymbol{x})}(\mathrm{d}h(\boldsymbol{x}))\right]\mathrm{d}u \ , \quad (60)$$

$$= \operatorname{tr}\left[\left(\int_0^1 e^{(1-u)h(\boldsymbol{x})}\boldsymbol{L}e^{uh(\boldsymbol{x})}\mathrm{d}u\right)(\mathrm{d}h(\boldsymbol{x}))\right] \ . \quad (61)$$

Choosing $h(\boldsymbol{x}) = -\Omega(\boldsymbol{x})$ and $\boldsymbol{L} = e^{\Omega(\boldsymbol{y})}$ leads to

$$\nabla_{\boldsymbol{x}}\operatorname{tr}(e^{\Omega(\boldsymbol{y})}e^{-\Omega(\boldsymbol{x})}) = -\Omega^*\int_0^1 e^{(u-1)\Omega(\boldsymbol{x})}e^{\Omega(\boldsymbol{y})}e^{-u\Omega(\boldsymbol{x})}\mathrm{d}u \ , \quad (62)$$

the integral term being Hermitian since it reads as the integral of Hermitian matrices

$$\int_0^{\frac{1}{2}} \underbrace{e^{(u-1)\Omega(\boldsymbol{x})}e^{\Omega(\boldsymbol{y})}e^{-u\Omega(\boldsymbol{x})} + e^{-u\Omega(\boldsymbol{x})}e^{\Omega(\boldsymbol{y})}e^{(u-1)\Omega(\boldsymbol{x})}}_{\text{Hermitian}}\mathrm{d}u \ .$$

## Acknowledgment

The authors would like to thank the anonymous reviewers as well as the associate editor for their comments, criticisms and encouragements. We also thank the Centre National d'Études Spatiales, the Office Nationale d'Études et de Recherches Aérospatiales and the Délégation Générale pour l'Armement

---

[10]http://www.math.u-bordeaux.fr/~cdeledal/mulog

| Image size | $D^2$ | File size | Gaussian denoiser | Time for GD | Time for NM | Total time |
|---|---|---|---|---|---|---|
| $256 \times 256$ | $2 \times 2$ | 2Mb | BM3D | 0.87s | 0.35s | 8.97s |
| | | | TV | 1.43s | | 17.08s |
| | | | DDID | 38.4s | | 4min 29s |
| $256 \times 256$ | $3 \times 3$ | 4.5Mb | BM3D | 0.87s | 2.20s | 29.21s |
| | | | TV | 1.43s | | 48.18s |
| | | | DDID | 38.4s | | 12min 6s |
| $512 \times 512$ | $3 \times 3$ | 18Mb | BM3D | 3.6s | 9.97s | 2min 09s |
| | | | TV | 6.3s | | 4min 4s |
| | | | DDID | 2min 37s | | 48min 52s |

for providing the RAMSES data, the Airbus Defence and Space for the TerraSAR-X data, the German Aerospace Center (DLR) for the F-SAR data.

This work has been carried out with financial support from the French State, managed by the French National Research Agency (ANR) in the frame of the GOTMI project (ANR-16-CE33-0010-01) and the financial support from DGA (Direction Générale de l'Armement) and ANR in the frame of ALYS project (ANR-15-ASTR-0002).

S. Tabti also thanks the DGA and the Centre National de la Recherche Scientifique for funding her PhD program.